# A study on Marangoni convection by the variational iteration method

**Onur Karaoğlu, Galip Oturanç**

Department of Mathematics, Faculty of Sciences, University of Selçuk, 42075, Turkey.

**Abstract**
In this letter, we will consider the use of the variational iteration method and Padé approximant for finding approximate solutions for a Marangoni convection induced flow over a free surface due to an imposed temperature gradient. The solutions are compared with the numerical (fourth-order Runge Kutta) solutions.

## 1. Introduction

Mass transfer across a fluid interface plays an important role in many chemical and engineering processes. Irregularities in mass transfer can lead to local changes in surface tension and density. Surface tension is a fundamental property of all liquids that causes the surface portion of a liquid to be attracted to another surface. Marangoni convection is the tendency of heat and mass in areas with low surface tension to move to areas with higher surface tension.

Marangoni convection has been investigated in other studies [1 – 5]. Slavtchev and Miladinova [1] presented similarity solutions for the thermocapillary flow in a thin layer of a liquid presenting the surface tension minimum. Arafune and Hirata [2] numerically studied interactive solutal and thermal Marangoni convections in a rectangular open boat. Okano et al. [3] presented a theoretical evaluation of natural and Marangoni convections using the finite difference method. Arafune and Hirata [4] measured the free surface velocity of thermal and solutal Marangoni convection in a rectangular open boat using In-Ga-Sb melt as the test fluid, Christopher and Wang [5] presented a similarity solution for Marangoni-induced flow both for the velocity profile and the temperature profile, assuming the development of boundary-layer flow along a surface with various temperature profiles imposed. Also, the subject has been investigated with a numerical method [6].

This paper presented a numerical analysis with a variational iteration method for Marangoni convection induced flow over a free surface due to an imposed temperature gradient.

## 2. Formulation of the problem

The laminar boundary layer equations for conservation of mass, momentum and energy for an incompressible viscous fluid are, respectively:

$$\frac{\partial u}{\partial x}+\frac{\partial \upsilon}{\partial y}=0 \qquad (1)$$

$$u\frac{\partial u}{\partial x}+\upsilon\frac{\partial u}{\partial y}=v\frac{\partial^2 u}{\partial y^2} \qquad (2)$$

$$u\frac{\partial T}{\partial x}+\upsilon\frac{\partial T}{\partial y}=\kappa\frac{\partial^2 T}{\partial y^2}. \qquad (3)$$

Here $v$ and $\kappa$ are the momentum and thermal diffusivities, respectively.

Nomenclature
- $u$      velocity component along x
- $\upsilon$      velocity component along y
- $\psi$      stream function
- $F$      dimensionless stream function
- $k$      parameter relating to the power law exponent

The boundary conditions at the free surface are:

$$\mu \left.\frac{\partial u}{\partial y}\right|_{y=0} = -\frac{d\sigma}{dT}\left.\frac{\partial T}{\partial x}\right|_{y=0} \tag{4a}$$

$$\upsilon(x,0) = 0 \tag{4b}$$

$$T(x,0) = T(0,0) + mx^{k+1}. \tag{4c}$$

Far from the surface, the velocity and thermal boundary conditions are:

$$u(x,\infty) = 0 \tag{5a}$$

$$\left.\frac{\partial T}{\partial y}\right|_{y=\infty} = 0, \tag{5b}$$

where $k = 0$ refers to a linear surface temperature profile, $k = 1$ is quadratic. The minimum value of $k$ is $-1$, which results in no temperature variation on the bubble surface and, thus, no Marangoni induced flow.

Using the stream function $\psi(x,y)$, the similarity variable $\eta$, the dimensionless stream function $F(\eta)$ and the dimensionless temperature function $\theta(\eta)$ can be introduced as

$$\eta = C_1 x^\alpha y \tag{6a}$$

$$F(\eta) = C_2 x^\beta \psi(x,y) \tag{6b}$$

$$\theta(\eta) = (T(x,y) - T(0,0)) x^t. \tag{6c}$$

The boundary layer equations (Eqs. (1) – (5)) with similarity transforms (6) can be combined into the following non-linear ordinary differential equations [5-7]

$$F''' = aF'^2 - bFF'' \tag{7a}$$

$$\theta'' = \Pr(-bF\theta' - tF'\theta) \tag{7b}$$

with the boundary conditions:

$$F(0) = 0, \quad F''(0) = -k-1, \quad F'(\infty) = 0 \tag{8}$$

$$\theta(0) = m, \quad \theta'(\infty) = 0, \tag{9}$$

where coefficients and parameters are:

$$C_1 = \sqrt[3]{\frac{(d\sigma/dT)m\rho}{\mu^2}}, \quad C_2 = \sqrt[3]{\frac{\rho^2}{(d\sigma/dT)m\mu}} \tag{10}$$

$$a = \frac{2k+1}{3}, \quad b = \frac{k+2}{3}, \quad t = -1-k. \tag{11}$$

Eqs. (7) - (8) are non-linear equations, and it is often difficult to find exact solutions for such equations. Therefore, the use of approximate solution methods for such equations becomes almost essential.

## 3. Variational iteration method

The variational iteration method was proposed by He [8–11].
To illustrate the basic concepts of the variational iteration method, consider the following differential equation:

$$L[u(\eta)] + N[u(\eta)] = g(\eta), \qquad (12)$$

where $L$ is a linear operator, $N$ is a non-linear operator and $g(x)$ is a given continuous function. According to the method, the correction functional can be constructed as follows,

$$u_{n+1}(\eta) = u_n(\eta) + \int_0^\eta \lambda(\tau)[Lu_n(\tau) + N\tilde{u}_n(\tau) - g(\tau)]d\tau, \qquad (13)$$

where $\lambda$ is a Lagrange multiplier [9,10,12], which can be optimally identified via variational theory [9 – 11], $u_n$ is the $n$th – order approximate solution, and $\tilde{u}_n$ is considered a restricted variation [9 – 11], i.e. $\delta \tilde{u}_n = 0$.

## 4. Solutions

In this section, Eqs. (7) – (9) are solved using the variational iteration method. First, let us start with solution for the momentum equation.
The correction functional for (7a) reads as

$$F_{n+1}(\eta) = F_n(\eta) + \int_0^\eta \lambda(\tau)\left(\frac{\partial^3 \tilde{F}_n}{\partial \tau^3} - a\left(\frac{\partial \tilde{F}_n}{\partial \tau}\right)^2 + b\tilde{F}_n\left(\frac{\partial^2 \tilde{F}_n}{\partial \tau^2}\right)\right)d\tau. \qquad (14)$$

The Lagrange multiplier can be identified as:

$$\lambda = -\frac{1}{2}(\tau - \eta)^2. \qquad (15)$$

Substituting Eq. (15) into Eq. (14) results in the following iteration formula:

$$F_{n+1}(\eta) = F_n(\eta) - \frac{1}{2}\int_0^\eta (\tau - \eta)^2 \left(\frac{\partial^3 F_n}{\partial \tau^3} - a\left(\frac{\partial F_n}{\partial \tau}\right)^2 + bF_n\left(\frac{\partial^2 F_n}{\partial \tau^2}\right)\right)d\tau. \qquad (16)$$

We start with an initial approximation with arbitrary unknown constants

$$F_0(\eta) = A + B\eta + Ce^{-\eta} \qquad (17)$$

where $A$, $B$ and $C$ are unknown constants to be further determined.
Substituting Eq. (17) into the iteration formula, Eq. (16), produces:

$$F_1(\eta) = B\eta + A + C - \frac{C^2}{24} - C\eta + \frac{\eta^2 C}{2} + \frac{2}{3}\eta e^{-\eta}CB - \frac{1}{12}\eta^2 C^2 - \frac{2}{3}CA - \frac{8}{3}BC + \frac{1}{12}\eta C^2 + \\ \frac{1}{18}B^2\eta^3 + \frac{1}{24}C^2 e^{-2\eta} + \frac{8}{3}BCe^{-\eta} + \frac{2}{3}CAe^{-\eta} + \frac{2}{3}CA\eta + 2CB\eta - \frac{2}{3}BC\eta^2 - \frac{1}{3}CA\eta^2 \qquad (18)$$

and so on.
Because of the difficulty in the calculation, this approach can be handled as a solution.
Consequently, using Eq. (18) yields the following approximations:

$$F'(\eta) = B - C + C\eta - \frac{2}{3}CB\eta e^{-\eta} - 2BCe^{-\eta} - \frac{1}{6}C^2 + \frac{1}{12}C^2 + \frac{1}{6}B^2\eta^2 - \frac{1}{12}C^2 e^{-2\eta} - \\ \frac{2}{3}CAe^{-\eta} + \frac{2}{3}CA + 2BC - \frac{4}{3}CB\eta - \frac{2}{3}CA\eta \qquad (19)$$

$$F''(\eta) = C + \frac{2}{3}CB\eta e^{-\eta} + \frac{4}{3}BCe^{-\eta} - \frac{1}{6}C^2 + \frac{1}{3}B^2\eta + \frac{1}{6}C^2 e^{-2\eta} + \frac{2}{3}CAe^{-\eta} - \frac{4}{3}BC - \frac{2}{3}CA \qquad (20)$$

Now, using the boundary conditions from Eq. (8) we will try to find the unknown coefficients in the approach of Eq. (18). However, the boundary conditions only give two of the three values that we need. We will try to find the other value through the condition of the $F'(\infty)=0$ by applying the Padé approximation to Eq. (19). The diagonal approximants in the technique of Padé are known to have more accurate approximations [7,13,14]. Here, [2/2] of Padé approximation is used.

When $k=0$, $a=1/3$, $b=2/3$, $t=-1$, we obtain $c=1.364053270$, and

$$F(\eta) = 1.437335891 + 0.1120407076\eta - 0.04691602733\eta^2 + 0.005155387494\eta^3 - 0.2030839727\eta e^{-\eta} - 1.479002557e^{-\eta} + \frac{1}{24}e^{-2\eta} \quad (21)$$

$$F'(\eta) = 0.1120407076 - 0.09383205466\eta + 0.01546616248\eta^2 - 0.2030839727\eta e^{-\eta} + 1.275918584e^{-\eta} - \frac{1}{12}e^{-2\eta} \quad (22)$$

$$F''(\eta) = -0.09383205466 + 0.03093232496\eta - 0.2030839727\eta e^{-\eta} - 1.072834611e^{-\eta} + \frac{1}{6}e^{-2\eta} \quad (23)$$

The standard fourth-order Runge-Kutta formulas and a shooting technique are used to test the reliability of the approximate solution. The behavior of the solutions are obtained by variational iteration method and the Runge-Kutta formulas are shown in Fig.1.

In order to solve Eq. (7b), we introduce the following definition
$$\theta(\eta) = mg(\eta)$$
so
$$g = \Pr(-bFg' - tF'g) \quad (24)$$
and the boundary condition in Eq. (9) reduces to:
$$g(0)=1, \quad g'(\infty)=0. \quad (25)$$

To solve (7b) by means of He's variational iteration method, a correction function can be written as follows:

$$g_{n+1}(\eta) = g_n(\eta) + \int_0^\eta \lambda(\tau)\left(\frac{\partial^2 g_n(\tau)}{\partial \tau^2} + \Pr bf(\tau)\left(\frac{\partial \tilde{g}_n(\tau)}{\partial \tau}\right)^2 + \Pr t\left(\frac{\partial F_n(\tau)}{\partial \tau}\right)g_n(\tau)\right)d\tau \quad (26)$$

where $\lambda$ is a Lagrange multiplier.

The Lagrange multiplier can be identified as:
$$\lambda = (\tau - \eta). \quad (27)$$

Substituting Eq. (27) into Eq. (26) results in the following iteration formula:

$$g_{n+1}(\eta) = g_n(\eta) + \int_0^\eta (\tau-\eta)\left(\frac{\partial^2 g_n(\tau)}{\partial \tau^2} + \Pr bf(\tau)\left(\frac{\partial \tilde{g}_n(\tau)}{\partial \tau}\right)^2 + \Pr t\left(\frac{\partial F_n(\tau)}{\partial \tau}\right)g_n(\tau)\right)d\tau. \quad (28)$$

Similarly, we start with an initial approximation with arbitrary unknown constants
$$g_0(\eta) = B\eta + Ce^{-\eta} \quad (29)$$

where $B$ and $C$ are unknown constants that need to be determined.

Comparison of the solutions of function $g'$ obtained using variational iteration method with [3/3] of Padé approximation and the Runge-Kutta formulas are shown in Fig. 2 for $\Pr=5$. As shown in the figures, this results in a solution to the limited ranges being discovered.

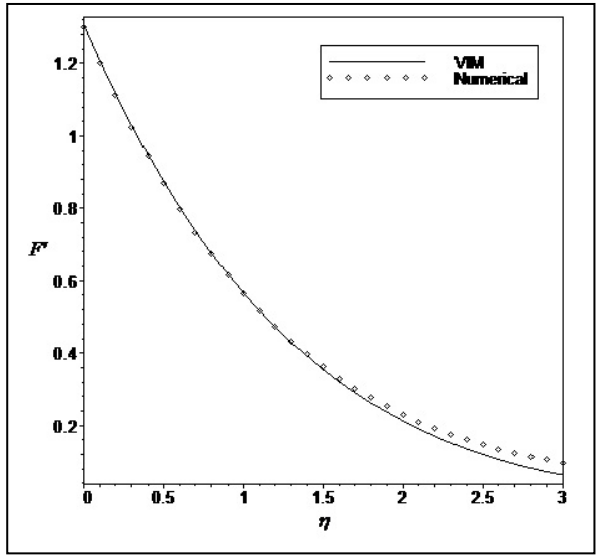 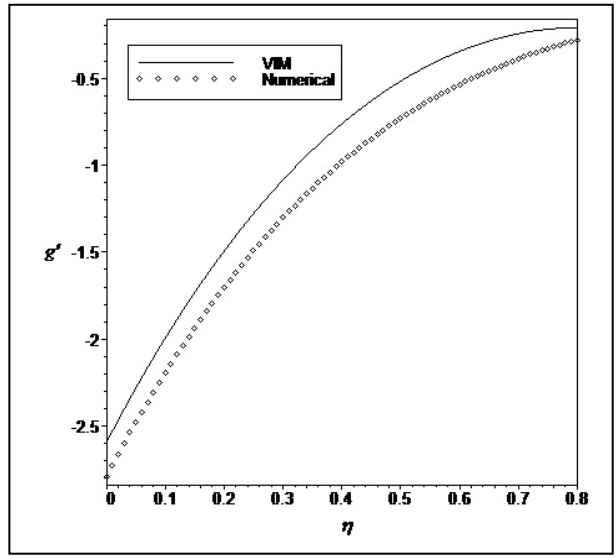

Fig.1. Velocity distribution profiles for $k = 0$.   Fig.2. Temperature gradient profiles for $k = 0$, $\Pr = 5$

## 5. Conclusions

In this paper, the variational iteration method and the Padé approximation technique were applied for solving Marangoni convection induced flow over a fluid gas-free surface cause by an imposed temperature gradient. The solutions were compared with the standard fourth-order Runge-Kutta solution and computations were performed using Maple13 package.

## Acknowledgements
This study is a part of Ph.D. Thesis the first author.